\documentclass[12pt]{amsart}
\usepackage{amsmath}
\usepackage{amssymb}
\usepackage{amscd}
\usepackage{latexsym}
\usepackage{enumerate}

\begin{document}

\numberwithin{equation}{section}
\renewcommand{\baselinestretch}{1.3}
\newtheorem{theorem}{Theorem}[section]
\newtheorem{proposition}[theorem]{Proposition}
\newtheorem{corollary}[theorem]{Corollary}
\newtheorem{lemma}[theorem]{Lemma}
\newtheorem{mlemma}[theorem]{Main Lemma}
\newtheorem{mtheorem}[theorem]{Main Theorem}
\newtheorem{blackbox}[theorem]{Black Box}
\newtheorem{example}[theorem]{Example}
\newtheorem{remark}[theorem]{Remark}
\newtheorem{rn}[theorem]{Remarks and Notations}
\newtheorem{construction}[theorem]{Construction}
\newtheorem{question}[theorem]{Question}
\newtheorem{definition}[theorem]{Definition}
\newtheorem{notation}[theorem]{Notation}
\newenvironment{prooff}{\noindent{\sc Proof.}}{\qed\medskip}
\newcommand{\Z}{{\mathbb Z}}
\newcommand{\Q}{{\mathbb Q}}
\newcommand{\N}{{\mathbb N}}
\newcommand{\K}{{\mathcal K}}
\newcommand{\Ka}{{\mathcal  K}}
\newcommand{\Ks}{{\mathcal K}^*}
\newcommand{\R}{{\mathbb R}}
\newcommand{\G}{{\mathcal G}}
\newcommand{\Pa}{{\mathcal P}}
\newcommand{\adic}{^{\raise0.05cm\hbox{\scriptsize $\wedge$}}}
\newcommand{\Kc}{{\mathcal K}\adic}
\newcommand{\wH}{\widehat H}
\newcommand{\wB}{\widehat B}
\newcommand{\wA}{\widehat A}
\newcommand{\wC}{\widehat C}
\newcommand{\Hom}{{\rm Hom}}
\newcommand{\Aut}{{\rm Aut}}
\newcommand{\Out}{{\rm Out}}
\newcommand{\End}{{\rm End}}
\newcommand{\Tor}{{\rm Tor}}
\newcommand{\Ext}{{\rm Ext}}
\newcommand{\Pext}{{\rm Pext}}
\newcommand{\coker}{{\rm coker}\,}
\newcommand{\im}{{\rm Im}\,}
\newcommand{\ident}{{\rm id}}
\newcommand{\dospunts}{\hspace{-0.1cm}:}
\newcommand{\red}{/\hspace{-0.05cm}/}
\newcommand{\Inn}{{\rm Inn}\,}
\newcommand{\pinn}{{\rm pInn}\,}
\newcommand{\dom}{{\rm dom}\,}
\newcommand{\cod}{{\rm cod}\,}
\newcommand{\Dom}{{\rm Dom}\,}
\newcommand{\cf}{{\rm cf}\,}
\newcommand{\br}{{\rm br}\,}
\newcommand{\Br}{{\rm Br}\,}
\newcommand{\infim}{{\rm inf}\,}
\newcommand{\PP}{{\mathbb P}}
\newcommand{\KK}{{\mathbb K}}
\newcommand{\GG}{{\mathbb G}}
\newcommand{\E}{{\mathbb E}}
\newcommand{\card}[1]{|#1|}
\newcommand{\norm}[1]{||#1||}
\newcommand{\Dl}{{\rm Dl}_\lambda}
\def\b{\beta}
\def\g{\gamma}
\def\d{\delta}
\def\e{\varepsilon}
\def\l{\lambda}
\def\lp{\lambda^+}
\def\blp{\widetilde{\lambda^+}}
\def\bu{\widetilde{u}}
\def\bv{\widetilde{v}}
\def\tx{\tilde{x}}
\def\k{\kappa}
\def\m{\mu}
\def\abar{\overline{\alpha}}
\def\aa{{\bf a}}
\def\ra{\rightarrow}
\def\arr{\longrightarrow}
\def\iff{\Longleftrightarrow}
\def\sigmaa{{\bf \Sigma_a}}
\def\L{{\mathfrak L}}
\def\F{{\mathfrak F}}
\def\c{{\mathfrak c}}
\def\X{{\mathfrak X}}
\def\P{{\mathfrak P}}
\def\p{{\mathfrak p}}
\def\Diam{\diamondsuit}
\def\mapdown#1{\Big\downarrow\rlap{$\vcenter{\hbox{$\scriptstyle#1}}$}}
\def\cRp{\rm{cRep}$_2R$ }
\def\cR{c$R_2$}
\def\restr{\upharpoonright}


\title{Almost-free $E$-rings of cardinality $\aleph_1$}

\author[R\"udiger G\"obel]{R\"udiger G\"obel}
\address{Fachbereich 6 -- Mathematik,
University of Essen, 45117 Essen, Germany}
\email{R.Goebel@uni-essen.de}
\author[Saharon Shelah]{Saharon Shelah}
\address{Institute of Mathematics, The Hebrew University, Givat Ram, Jerusalem 91904,
Israel} \email{shelah@math.huji.ac.il}
\author[Lutz Str\"ungmann]{Lutz Str\"ungmann}
\address{Fachbereich 6 -- Mathematik,
University of Essen, 45117 Essen, Germany}

\curraddr{Institute of Mathematics, The Hebrew University, Givat
Ram, Jerusalem 91904, Israel} \email{lutz@math.huji.ac.il}
\thanks{
Publication 785 in the second authors list of publication. The
third author was supported by a MINERVA fellowship.}


\begin{abstract}
An $E$-ring is a unital ring $R$ such that every endomorphism of
the underlying abelian group $R^+$ is multiplication by some
ring-element. The existence of almost-free $E$-rings of
cardinality greater than $2^{\aleph_0}$ is undecidable in ZFC.
While they exist in Goedel's universe, they do not exist in other
models of set theory. For a regular cardinal $\aleph_1 \leq
\lambda \leq 2^{\aleph_0}$ we construct $E$-rings of cardinality
$\lambda$ in ZFC which have $\aleph_1$-free additive structure.
For $\lambda=\aleph_1$ we therefore obtain the existence of
almost-free $E$-rings of cardinality $\aleph_1$ in ZFC.
\end{abstract}

\maketitle


\section{Introduction}
Recall that a unital ring $R$ is an $E$-ring if the evaluation map
$\varepsilon :$ $\End_{\Z}(R^+) \rightarrow R$ given by $\varphi
\mapsto \varphi(1)$ is a bijection. Thus every endomorphism of the
abelian group $R^+$ is multiplication by some element $r \in R$.
$E$-rings were introduced by Schultz \cite{S} and easy examples
are subrings of the rationals $\Q$ or pure subrings of the ring of
$p$-adic integers. Schultz characterized $E$-rings of finite rank
and the books by Feigelstock \cite{Feig1}, \cite{Feig2} and an
article \cite{PV} survey the results obtained in the eighties, see
also \cite{Re}, \cite{F}. In a natural way the notion of $E$-rings
extends to modules by calling a left $R$-module $M$ an
$E(R)$-module or just $E$-module if $\Hom_{\mathbb{Z}}(R,M) =$
$\Hom_R(R,M)$ holds, see \cite{BS}. It turned out that a unital
ring $R$ is an $E$-ring if and only if it is an $E$-module.\\
$E$-rings and $E$-modules have played an important role in the
theory of torsion-free abelian groups of finite rank. For example
Niedzwecki and Reid \cite{NR} proved that a torsion-free abelian
group $G$ of finite rank is cyclically projective over its
endomorphism ring if and only if $G=R \oplus A$, where $R$ is an
$E$-ring and $A$ is an $E(R)$-module. Moreover, Cassacuberta and
Rodr\'{\i}guez \cite{CRT}
noticed the role of $E$-rings in homotopy theory.\\
It can be easily seen that every $E$-ring has to be commutative
and hence can not be free as an abelian group except when $R =\Z$.
But it was proved in \cite{DMV}, using a Black Box argument, that
there exist arbitrarily large $E$-rings $R$ which are
$\aleph_1$-free as a group which means that every countable
subgroup of $R^+$ is free. This implies the existence of
$\aleph_1$-free $E$-rings of cardinality $\aleph_1$ under the
assumption of the continuum hypothesis. Moreover, it was shown in
\cite{GSt} that there exist almost-free $E$-rings for any regular
not weakly compact cardinal $\kappa > \aleph_0$ assuming diamond,
a prediction principle which holds for example in Goedel's
constructible universe. Here, a group of cardinality $\lambda$ is
called almost-free if all its subgroups of smaller cardinality
than $\lambda$ are free.\\ Since the existence of $\aleph_2$-free
$E$-rings of cardinality $\aleph_2$ is undecidable in ordinary set
theory ZFC (see \cite[Theorem 5.1]{GS2} and \cite{GSt}) it is
hopeless to conjecture that there exist almost-free $E$-rings of
cardinality $\kappa$ in ZFC for cardinals $\kappa$ larger than
$2^{\aleph_0}$. However, we will prove in this paper that there
are $\aleph_1$-free $E$-rings in ZFC of cardinality $\lambda$ for
every regular cardinal $\aleph_1 \leq \lambda \leq 2^{\aleph_0}$.
Thus we show the existence of almost-free $E$-rings of size
$\aleph_1$ in ZFC. \\
The construction of $\aleph_1$-free $E$-rings $R$ of cardinality
$\aleph_1$ in ZFC is much easier if $|R| = 2^{\aleph_0}$, this
because in the first case we are closer to freeness which tries to
prevent endomorphisms to be scalar multiplication. The underlying
setting and the combinatorial predictions however are similar to
the one used by the first two authors in \cite{GS1} for
constructing indecomposable almost-free groups of cardinality
$\aleph_1$ with prescribed endomorphism ring.\\

Our notations are standard and for unexplained notions we refer to
\cite{F1, F2, F3} for abelian group theory and to \cite{EM} for
set-theory. All groups under consideration are abelian.

\section{Topology, trees and a forest}
In this section we explain the underlying geometry of our
construction which was used also in \cite{GS1}, see  \cite{GS1}
for further details.\\

Let $F$ be a fixed countable principal ideal domain with a fixed
infinite set $S=\{ s_n : n \in \omega \}$ of pair-wise coprime
elements, that is $s_nF + s_mF=F$ for all $n\ne m$, for brevity
say $F$ is a {\em $p$-domain}. We chose a sequence of elements
\[ q_0 =1 \text{ and } q_{n+1}= s_nq_n \text{ for all } n \in \omega
\]
in $F$, hence the descending chain $q_nF \ (n\in \omega)$ of
principal ideals satisfies $\bigcap\limits_{n \in \omega} q_nF =0$
and generates the Hausdorff $S$-topology on $F$.\\

Now let $T=$ $^{\omega>}2$ denote the tree of all finite branches
$\tau : n \longrightarrow 2$ $(n \in \omega)$. Moreover,
$^{\omega}2= \Br(T)$ denotes all infinite branches $\eta: \omega
\longrightarrow 2$ and clearly $\eta\restr_n \in T$ for all $\eta
\in  \Br(T)$ $(n \in \omega)$. If $\eta \not= \mu \in \Br(T)$ then
\[ \br(\eta,\mu)= \infim \{ n \in \omega : \eta(n) \not= \mu(n) \} \]
denotes the {\it branch point} of $\eta$ and $\mu$. If $C \subset
\omega$ then we collect the subtree
\[ T_C = \{ \tau \in T : \text{ if } e \in l(\tau) \backslash C \text{ then }
\tau(e)=0 \} \] of $T$ where $l(\tau)=n$ denotes the {\it length} of the
finite branch $\tau:n
\longrightarrow 2$.\\
Similarly,
\[ \Br(T_C) = \{ \eta \in \Br(T) : \text{ if } e \in \omega
\backslash C \text{ then } \eta(e)=0 \} \] and hence $\eta
\restr_n \in T_C$ for all $\eta \in \Br(T_C)$ $(n \in
\omega)$.\\
Now we collect some trees to build a forest. Therefore, let $\aleph_1 \leq \lambda \leq 2^{\aleph_0}$ be a regular cardinal.\\
Choose a family $\mathfrak{C} = \{ C_{\alpha} \subset \omega :
\alpha < \lambda \}$ of pair-wise almost disjoint infinite subsets
of $\omega$. Let $T \times \alpha = \{ v \times \alpha : v \in T
\}$ be a disjoint copy of the tree $T$ and let $T_{\alpha} =
T_{C_{\alpha}} \times \alpha$ for $\alpha < \lambda$. For
simplicity we denote the elements of $T_{\alpha}$ by $\tau$
instead of $\tau \times \alpha$ since it will always be clear from
the context to which $\alpha$ the finite branch $\tau$ refers to.
By \cite[Observation 2.1]{GS1} we may assume that each tree
$T_{\alpha}$ is perfect for $\alpha < \lambda$, i.e. if $n \in
\omega$ then there is at most one finite branch $\eta\restr_n$
such that $\eta\restr_{(n+1)}\ \not= \mu\restr_{(n+1)}$ for some
$\mu \in T_{\alpha}$. We build a forest by letting
\[ T_\Lambda=\bigcup\limits_{\alpha < \lambda}T_{\alpha}. \] Now
we define our {\it base algebra} as $B_\Lambda=F[z_{\tau} : \tau
\in T_\Lambda ]$ which is a pure and dense subalgebra of its
$S$-adic completion $\widehat{B_\Lambda}$ taken in the
$S$-topology on $B_\Lambda$. \\
For later use we state the following definition which allows us to
view the algebra $B_\Lambda$ as a module generated over $F$ by
monomials in the ``variables'' $z_{\tau}$ ($\tau \in T_\Lambda$).

\begin{definition}
Let $X$ be a set of commuting variables and $R$ an $F$-algebra.
Then any map $\sigma: X \rightarrow R$ extends to a unique
epimorphism (also called) $\sigma: F[X] \rightarrow F[\sigma(X)]$.
Thus any $r \in F[\sigma(X)]$ can be expressed by a polynomial
$\sigma_r \in F[X]$, which is a preimage under $\sigma$: There are
$l_1,\cdots,l_n$ in $\sigma(X)$ such that $ r = \sigma_r(l_1,
\cdots, l_n)$ becomes a polynomial-like expression. Similarly
$M(\sigma(X))$ is the monomial-like set of all products taken from
$\sigma(X)$.
\end{definition}

In particular, if $Z_{\alpha}=\{ z_{\tau} : \tau \in T_{\alpha}
\}$ ($\alpha < \lambda$) and $Z_\Lambda=\{ z_{\tau} : \tau \in
T_\Lambda \}$, then as always the polynomial ring $B_\Lambda$ can
be viewed as a free $F$-module over the basis of monomial, we
have $B_\Lambda=\bigoplus\limits_{z \in M(Z_\Lambda)}zF$.\\
Since $\aleph_1 \leq \lambda \leq 2^{\aleph_0}= |\Br(T_{C_\alpha
})|$ we can choose a family $\{V_\alpha  \subseteq$
$\Br(T_{C_\alpha }) : \alpha < \lambda \}$ of subsets $V_\alpha $
of $\Br(T_{C_\alpha })$ with $|V_\alpha |=\lambda$ for ($\alpha <
\lambda$). Note that for $\alpha \ne \beta < \lambda$ the infinite
branches from $V_\alpha $ and $V_{\beta}$ branch at almost
disjoint sets since $C_\alpha \cap C_{\beta}$ is finite, thus the
pairs $V_\alpha $, $V_{\beta}$ are disjoint. Moreover, we may
assume that for any $m \in \omega$, $\lambda$ pairs of branches in
$V_\alpha $ branch above $m$.

\section{The Construction}

Following \cite{GS1} we use the definition\\

\begin{definition}
Let $x \in \widehat{B_\Lambda}$ be any element in the completion
of the base algebra $B_\Lambda$. Moreover, let $\eta \in V_\alpha
$ for $\alpha < \lambda$. Then we define the following elements
for $n \in \omega$:
\[ y_{\eta nx}:= \sum\limits_{i \geq n} \frac{q_i}{q_n} \left(
z_{\eta\restr_i} \right) + x\sum\limits_{i \geq n} \frac{q_i}{q_n}
\eta(i). \] The elements $y_{\eta nx}$ are called {\it branch like
elements}.
\end{definition}

Note that each element $y_{\eta nx}$ connects an infinite branch
$\eta \in \Br(T_{C_\alpha })$ with finite branches from the
disjoint tree $T_\alpha $. Furthermore, the element $y_{\eta nx}$
encodes the infinite branch $\eta$ into an element of
$\widehat{B_\Lambda}$. We have a first observation which describes
this as an equation and which is crucial for the rest of this
paper.

\begin{eqnarray} \label{equat} y_{\eta nx}=s_{n+1} y_{\eta (n+1) x} +
z_{\eta\restr_n} + x\eta(n) \mbox{ for all }\alpha < \lambda,
\eta \in V_\alpha . \end{eqnarray}

\proof We calculate the difference $q_ny_{\eta n x} -
q_{n+1}y_{\eta (n+1)x}$:
\begin{alignat}{2}
q_ny_{\eta n x} - q_{n+1}y_{\eta (n+1)x} &= \sum\limits_{i \geq n}
q_i \left( z_{\eta\restr_i} \right) + x\sum\limits_{i \geq n} q_i
\eta(i) - \sum\limits_{i \geq n+1} q_i \left( z_{\eta\restr_i}
\right) - x\sum\limits_{i \geq
n+1} q_i \eta(i) \\
& = q_nz_{\eta\restr_n} + q_nx\eta(n). \notag
\end{alignat}
Dividing by $q_n$ yields $y_{\eta nx}=s_{n+1} y_{\eta (n+1) x} +
z_{\eta\restr_n} + x\eta(n).$ \qed

The elements of the polynomial ring $B_\Lambda$ are unique finite
sums of monomials in $Z_\lambda$ with coefficients in $F$. Thus,
by $S$-adic topology, any $0\ne g \in \widehat{B_\Lambda}$ can be
expressed uniquely as a sum
\[ g= \sum\limits_{z \in [g]} g_{z}, \]
where $z$ runs over an at most countable subset $[g]\subseteq
M(Z_\Lambda)$ of monomials and $0 \ne g_z\in z\widehat{F}$.  We
put $[g]= \emptyset$ if $g = 0$. Thus any $g\in
\widehat{B_\Lambda}$ has a unique {\em support} $[g] \subseteq
M(Z_\Lambda)$, and support extends naturally to subsets of
$\widehat{B_\Lambda}$ by taking unions of the support of its
elements.  It follows that
\[ [y_{\eta no}]=\{ z_{\eta\restr_j \times \alpha} : j \in \omega,
j\geq n \} \] for any $\eta\in V_\alpha, \ n\in \omega$ and $[z] =
\{z\}$ for any $z\in M(Z_\Lambda)$.

Support can be used to define the norm of elements. If $X
\subseteq M(Z_\Lambda)$ then $$||X||=\infim\{ \beta < \lambda : X
\subseteq \bigcup\limits_{\alpha < \beta} M(Z_\alpha ) \}$$ is the
{\it norm} of $X$. If the infimum is taken over an unbounded
subset of $\lambda$, we write $||X||=\infty$. However, since
$\cf(\lambda) >\omega$, the {\em norm of an element} $g \in
B_\Lambda$ is $||g||=||[g]||< \infty$ which is an ordinal $<
\lambda$ hence either discrete or cofinal to $\omega$. Norms
extend naturally to subsets of  $B_\Lambda$. In particular
$||y_{\eta
no}||=\alpha +1$ for any $\eta \in V_\alpha $.\\

We are ready to define the final $F$-algebra $R$ as a
$F$-subalgebra of the completion of $B_\Lambda$. Therefore choose
a transfinite sequence $b_\alpha $ $(\alpha < \lambda)$ which runs
$\lambda$ times through the non-zero pure elements $b \in
B_\Lambda$ which are of the form $b=\sum\limits_{m \in M}m$ where
$M$ is a finite subset of $M(T_\Lambda)$. Note that $B_\Lambda/Fb$
is a free $F$-module.

\begin{definition}
Let $F$ be a countable principal ideal $p$-domain with identity
and let $B_\Lambda:= F [ z_{\tau} : \tau \in T_\Lambda ]$ be the
polynomial ring over $Z_\Lambda$ as above. Then we define the
following smooth ascending chain of $F$-subalgebras of
$\widehat{B_\Lambda}$:

\begin{enumerate}

\item $R_0= \{0\}$; $R_1 :=F$;

\item $R_\alpha = \bigcup\limits_{\beta < \alpha} R_{\beta}$, for $\alpha$
a limit ordinal;

\item $R_{\alpha + 1} = R_\alpha  [ y_{\eta
nx_\alpha }, z_{\tau} :  \eta \in V_\alpha , \tau \in T_\alpha , n
\in \omega ] $;

\item $R = R_{\lambda} = \bigcup\limits_{\alpha < \lambda} R_\alpha $.

\end{enumerate}

We let $x_\alpha =b_\alpha $ if $b_\alpha  \in R_\alpha $ with
$||b_\alpha || \leq \alpha$ and $x_\alpha =0$ otherwise.
\end{definition}

For the rest of this paper purification and properties like
freeness, linear dependence or rank are taken with respect to $F$.
First we prove some properties of the ring $R_\alpha $ ($\alpha
\le \lambda$). It is easy to see that $R_\alpha =F[y_{\eta n
x_\alpha }, z_{\tau} : \eta \in V_{\beta}, \tau \in T_{\beta}, n
\in \omega, \beta < \alpha]$ is not a polynomial ring.
Nevertheless we have the following

\begin{lemma} \label{linunabh} For any fixed $n\in\omega$ and
$\alpha < \lambda$ the set $M \left(y_{\eta nx_\alpha }, z_{\tau}
: \eta \in V_\alpha , \tau \in T_\alpha  \right)$ is linearly
independent over $R_\alpha$. Thus $R_\alpha [y_{\eta n x_\alpha },
z_{\tau}: \eta \in V_\alpha , \tau \in T_\alpha ]$ is a polynomial
ring.
\end{lemma}

\proof Assume that the set $M \left( y_{\eta nx_\alpha }, z_{\tau}
: \eta \in V_\alpha , \tau \in T_\alpha  \right)$ is linearly
dependent over $R_\alpha $ for some $\alpha < \lambda$ and $n \in
\omega$. Then there exists a non-trivial linear combination of the
form:\begin{equation} \label{*}
 \sum\limits_{y \in Y} \sum\limits_{z \in E_y} g_{y,z} y z = 0
\end{equation}
with $g_{y,z} \in R_\alpha $ and finite sets $Y \subset M \left(
y_{\eta nx_\alpha } : \eta \in V_\alpha  \right)$ and $E_y \subset
M \left( Z_\alpha  \right)$. We have chosen $V_{\beta} \cap
V_{\gamma} = \emptyset$ for all $\beta \not= \gamma$ and
$M(Z_\alpha ) \cap R_\alpha  = \emptyset$. Moreover $\| R_\alpha
\| < \| R_{\alpha +1} \|$ and hence there exists a basal element
$z_y \in B_\Lambda$ for any $1 \not= y \in Y$ with the following
properties:
\renewcommand{\theenumi}{\roman{enumi}}\renewcommand{\labelenumi}
{(\theenumi)}\begin{enumerate}\item $z_y \not\in E_{\tilde{y}}$
for all $\tilde{y} \in Y$;\item $z_y \not\in [\tilde{y} ]$ for all
$y \not= \tilde{y} \in Y$;\item $z_y \not\in [g_{\tilde{y},z}]$
for all $\tilde{y} \in Y, z \in E_{\tilde{y}}$;\item $z_y \in
[y]$.\end{enumerate}

Now we restrict the equation (\ref{*}) to the basal element $z_y$
and obtain $g_{y,z} z_y z=0$ for all $z \in E_y$. Since $z_y
\not\in [g_{y, z}]$ we derive $g_{y,z} =0$ for all $1\not=y \in Y$
and $z \in E_y$. Therefore equation (\ref{*}) reduces to
$\sum\limits_{z \in E_1} g_{1,z} z = 0$. We apply $M(Z_\alpha )
\cap R_\alpha  = \emptyset$ once more. Since each $z$ is a basal
element from the set $M(Z_\alpha )$ we get, equating coefficients,
that $g_{1,z}=0$ for all $z \in E_1$. Hence $g_{y,z}=0$ for all $y
\in Y, z \in E_y$, contradicting the assumption that (\ref{*}) is
a non-trivial linear combination. \qed

The following lemma shows that the $F$-algebras
$R_{\delta}/s_{n+1}R_{\delta}$ are also polynomial rings over
$F/s_{n+1}F$ for every $n < \omega$. Therefore, choose for $\delta
< \lambda$ and $V \subseteq V_{\delta}$, $|V|=\lambda$ a minimal
countable subset $U_V \subseteq V$ such that
\begin{equation}
\label{representatives}
 (\forall \eta \in V)(\forall n \in \omega)(\exists \rho \in U_V)
\text{ such that } \eta\restr_n \ = \rho\restr_n.
\end{equation}
In the next lemma we consider the particular case $V=V_{\delta}$
and we let $U_{V_{\delta}}=U_{\delta}$. The more general $U_V$
will be used later on.
\begin{lemma}
\label{polyring} Let $n < \omega$. Then the quotient algebra
$R_{\delta}/s_{n+1}R_{\delta}$ over $F/s_{n+1}F$ is generated by
the set of linearly independent elements
\[ X^{\delta}_{n+1}=\{ y_{\eta n x_{\beta}}, y_{\eta (n+1)x_{\beta}},
z_{\tau} : \eta \in U_{\beta}, \tau \in T_{\beta}, \tau \not=
\eta\restr_{n} \beta < \lambda \}. \] Thus
$R_{\delta}/s_{n+1}R_{\delta} = F/s_{n+1}F[y_{\eta n x_{\beta}},
y_{\eta (n+1)x_{\beta}}, z_{\tau} : \eta \in U_{\beta}, \tau \in
T_{\beta}, \tau \not= \eta\restr_{n} \beta < \lambda]$ is a
polynomial ring.
\end{lemma}

\proof The proof that the elements of $M(X^{\delta}_{n+1})$ are
linearly independent over $F/s_{n+1}F$ uses similar arguments as
in the proof of Lemma \ref{linunabh} and we therefore omit it.
Note that the only dependence relations between $y_{\eta n
x_{\beta}}$ and $y_{\eta (n+1)x_{\beta}}$ come from the equation
(\ref{equat}) and are therefore avoided by the definition of the
set $X^{\delta}_{n+1}$. Thus we only have to prove that
$R_{\delta}/s_{n+1}R_{\delta}=\left(F/s_{n+1}F\right)[X^{\delta}_{n+1}]$.
We will show by induction on $\alpha < \delta$ that
\[ \left(R_\alpha  + s_{n+1}R_{\delta} \right)/s_{n+1}R_{\delta}
\subseteq \left(F/s_{n+1}F \right)[X^{\delta}_{n+1}]. \] If
$\alpha=0$ or $\alpha =1$ then the claim is trivial, hence assume
that $\alpha >1$ and for all $\beta < \alpha$ we have \[
\left(R_{\beta} + s_{n+1}R_{\delta} \right)/s_{n+1}R_{\delta}
\subseteq \left(F/s_{n+1}F \right)[X^{\delta}_{n+1}]. \]

If $\alpha$ is a limit ordinal, then $\left(R_\alpha  +
s_{n+1}R_{\delta} \right)/s_{n+1}R_{\delta} \subseteq
\left(F/s_{n+1}F \right)[X^{\delta}_{n+1}]$ is immediate. Thus
assume that $\alpha = \beta +1$. By assumption and $x_{\beta} \in
R_{\beta}$ we know that $\left( x_{\beta} + s_{n+1}R_{\delta}
\right) \in \left(F/s_{n+1}F \right)[X^{\delta}_{n+1}]$. Hence
equation (\ref{equat}) shows that the missing elements
$z_{\eta\restr_{n}} + s_{n+1}R_{\delta}$ $(\eta \in V_{\beta})$
are in $\left( F/s_{n+1}F \right)[X^{\delta}_{n+1}]$. By induction
on $m< \omega$ using (\ref{equat}) it is now easy to show that
also $y_{\eta mx_{\beta}} + s_{n+1}R_{\delta} \in \left(
F/s_{n+1}F \right)[X^{\delta}_{n+1}]$ for every $m < \omega$ and
$\eta \in V_{\beta}$ and hence $R_\alpha  + s_{n+1}R_{\delta}
\subseteq \left( F/s_{n+1}F \right)[X^{\delta}_{n+1}]$ which
finishes the proof. \qed

Now we are able to prove that the members $R_\alpha $ of the chain
$\{ R_{\sigma} : \sigma < \lambda \}$ are $F$-pure submodules of
 $R$ and that $R$ is an $\aleph_1$-free domain.\\

\begin{lemma}
\label{domain} $R$ is a commutative $F$-algebra without
zero-divisors and $R_\alpha $ is $F$-pure in $R$ for all $\alpha <
\lambda$.
\end{lemma}

\proof By definition each $R_\alpha $ is a commutative $F$-algebra
and hence $R$ is commutative. To show that $R$ has no
zero-divisors it is enough to show that each member $R_\alpha $ of
the chain $\{ R_{\sigma} : \sigma < \lambda \}$ is an $F$-algebra
without zero-divisors. Since $F$ is a domain we can assume, by
induction, that $R_{\beta}$ has no zero-divisors for all $\beta<
\alpha$ and some $1 < \alpha < \lambda$. If $\alpha$ is a limit
ordinal then it is immediate that $R_\alpha $ has no
zero-divisors. Hence $\alpha = \gamma +1$ is a successor ordinal
and $R_{\gamma}$ is a domain. If $g,h \in R_\alpha $ with $gh = 0
\ne g $, then we must show that $h=0$. Write $g$ in the form
\begin{equation}
 g = \sum\limits_{y \in Y_g} \sum\limits_{z \in E_{g,y}} g_{y,z} y z
 \tag{g}
\end{equation}
with $0 \not= g_{y,z} \in R_{\gamma}$ and finite sets $Y_g \subset
M \left( y_{\eta nx_{\gamma}} : \eta \in V_{\gamma} \right)$ for
some fixed $n \in \omega$ and $E_{g,y} \subset M \left( Z_{\gamma}
\right)$. By (\ref{equat}) and $x_{\gamma} \in R_{\gamma}$ we may
assume $n$ is fixed. Similarly, we write
\begin{equation}
 h= \sum\limits_{y \in Y_h} \sum\limits_{z \in E_{h,y}} h_{y,z} y z
 \tag{h}
\end{equation}
with $h_{y,z} \in R_{\gamma}$ and finite sets $Y_h \subset M
\left( y_{\eta nx_{\gamma}} : \eta \in V_{\gamma} \right)$ and
$E_{h,y}
\subset M \left( Z_{\gamma} \right)$.\\
Next we want $h_{y,z} =0$ for all $y \in Y_h, z \in E_{h,y}$. The
proof follows by induction on the number of $h_{y,z}$'s. If $h=
h_{w,z^{\prime}}wz^{\prime}$, then
\[ gh= \sum\limits_{y \in Y_g,  z \in
E_{g,y}} g_{y,z}h_{w,z^{\prime}} y z w z^{\prime} \] and from
Lemma \ref{linunabh} follows $g_{y,z}h_{w,z^{\prime}} = 0$ for all
$y \in Y_g, z \in E_{g,y}$. Since $R_{\gamma}$ has no
zero-divisors we obtain $h_{w,z^{\prime}}=0$ and thus $h=0$. Now
assume that $k+1$ coefficients $h_{y,z}\ne 0$ appear in (h). We
fix an arbitrary coefficient $h_{w,z^{\prime}}$ and write $h=
h_{w,z^{\prime}} w z^{\prime} + h^{\prime}$ so that $wz'$ does not
appear in the representation of $h^{\prime}$. Therefore the
product $gh$ is of the form
\begin{equation}
 gh= \sum\limits_{y \in Y_g} \sum\limits_{z \in E_{g,y}}
g_{y,z}h_{w,z^{\prime}} y z w z^{\prime}  + gh^{\prime}. \tag{gh}
\end{equation}
If the monomial $wz^{\prime}$ appears in the representation of (g)
then the monomial $w^2(z^{\prime})^2$ appears in the
representation of (gh) only once with coefficient
$g_{w,z^{\prime}}h_{w,z^{\prime}}$. Using Lemma \ref{linunabh} and
the hypothesis that $R_{\gamma}$ has no
zero-divisors we get $h_{w,z^{\prime}}=0$.\\
If the monomial $w z^{\prime} $ does not appear in the
representation of (g) then $g_{y,z}h_{w,z^{\prime}}=0$ for all
appearing coefficients $g_{y,z}$ is immediate by Lemma
\ref{linunabh}. Thus $h_{w,z^{\prime}}=0$ and $h=h^{\prime}$
follows. By induction hypothesis also $h=0$ and $R$ has no zero-divisors.\\

It remains to show $F$-purity of $R_\alpha $ in $R$ for $\alpha <
\lambda$. Let $g \in R \setminus R_\alpha $ such that $fg \in
R_\alpha $ for some $0\ne f \in F$. If $\beta < \lambda$ is
minimal with $g \in R_{\beta}$ then $\beta > \alpha$ and it is
immediate that $\beta = \gamma +1$ for some $\gamma \ge \alpha$,
hence $fg \in R_\alpha \subset R_{\gamma}$. No we can write
\begin{equation}g = \sum\limits_{y \in Y_{g}} \sum\limits_{z \in E_{g,y}}
g_{y,z} y z \tag{g}
\end{equation}
with $g_{y,z} \in R_{\gamma}$ and finite sets $Y_{g} \subset M
\left( y_{vkx_{\gamma}} : v \in V_{\gamma} \right)$ for some fixed
$k \in \omega$ and $E_{g} \subset M \left( Z_{\gamma} \right)$ and
clearly
\[ fg=  \sum\limits_{y \in Y_{g}} \sum\limits_{z \in E_{g,y}}
fg_{y,z} y z \in R_{\gamma}. \] Hence there exists
$g_{\gamma} \in R_{\gamma}$ such that
\[ fg - g_{\gamma}=  \sum\limits_{y \in Y_{g}} \sum\limits_{z \in
E_{g,y}} fg_{y,z} y z - g_{\gamma} = 0. \] From Lemma
\ref{linunabh} follows $fg_{y,z}=0$ for all $ 1\not=y\in Y_g,
1\not=z \in E_{g,y}$, and $g_{y,z}=0$ because $R$ is
$F$-torsion-free. Hence (g) reduces to the summand with $y=z=1$,
but $g= g_{1,1} \in R_{\gamma}$ contradicts the minimality of
$\beta$. Thus $g \in R_\alpha $ and $R_\alpha $ is pure in
$R$.\qed

\begin{theorem}
\label{aleph1frei} Let $F$ be a countable principal ideal
$p$-domain with identity. If $R=\bigcup\limits_{\alpha < \lambda}
R_\alpha $ is the $F$-algebra constructed above then $R$ is a
domain such that $|R| = \lambda$ and $R$ is an $\aleph_1$-free
$F$-module. Moreover, for every ordinal $\alpha < \lambda$, the
quotient $R/R_\alpha $ is $\aleph_1$-free.\end{theorem}

\proof $| R | = \lambda$ is immediate by construction and $R$ is a
domain by Lemma \ref{domain}. It remains to show that $R$ is an
$\aleph_1$-free $F$-module. By Pontryagin's Theorem (see \cite[p.
93, Theorem 19.1]{F1}) it is enough to show that any pure
submodule of finite rank is contained in a free submodule.
Therefore let $U \subseteq R$ be a pure submodule of finite rank.
There exist elements $u_i \in R$ such that
\[ U= \left< u_1,...,u_n \right>_* \subseteq R. \]
Hence there is a minimal $\alpha < \lambda$ such that $u_i \in
R_\alpha $ for $i =1,...,n$, which obviously is a successor
ordinal $\alpha = \gamma +1$. Moreover, $U \subseteq R_\alpha $
since $R_\alpha $ is pure in $R$ and by induction we may assume
that $R_{\gamma}$ is $\aleph_1$-free. Using $ R_\alpha  =
R_{\gamma +1} = R_{\gamma}[ y_{\eta mx_{\gamma}}, z_{\tau} : \eta
\in V_{\gamma}, \tau \in T_{\gamma} m \in \omega] $ we can write
\[ u_i = \sum\limits_{y \in Y_i} \sum\limits_{z \in E_{i,y}} g_{y,z,i}
y z  \] with $g_{y,z,i} \in R_{\gamma}$ and finite sets $Y_i
\subset M \left( y_{\eta mx_{\gamma}} : \eta \in V_{\gamma}
\right)$ for some fixed $m
\in \omega$ and $E_{i,y} \subset M \left( Z_{\gamma} \right)$.\\
Choose the pure submodule
\[ R_{U} := \left< g_{y,z,i} : y \in Y_i, z \in E_{i,y}, 1 \leq i \leq
n \right>_* \subseteq R_{\gamma}\] of $R_{\gamma}$ and build the
set
\[ U^{\prime} := \left\{ y, z : y \in Y_i, z \in E_{i,y}, 1 \leq i
\leq n \right\}. \] Hence $\left< U^{\prime} \right>_{R_U}
\subseteq_* R_\alpha $ by Lemma \ref{linunabh} and purity of $R_U$
in $R_{\gamma}$. Thus $U \subseteq_* \left< U^{\prime}
\right>_{R_U} \subseteq_* R_\alpha $ and it remains to show that
$\left< U^{\prime} \right>_{R_U}$ is a free $F$-module. By
assumption $R_{\gamma}$ is $\aleph_1$-free and $R_U$ is a pure
submodule of finite rank of $R_{\gamma}$, hence $R_U$ is free. By
Lemma \ref{linunabh} we know that $U^{\prime}$ is linearly
independent over $R_{\gamma}$ and thus also over $R_U$. Now
$F$-freeness of $\left< U^{\prime} \right>_{R_U}$ follows and the
proof is complete. Similarly arguments show that for every $\alpha
< \lambda$ the quotient $R/R_\alpha $ is $\aleph_1$-free.\qed

\section{Main Theorem} In this final section we will prove
that the constructed $F$-algebra $R$ from Section $3$ is an
$E(F)$-algebra, i.e. that every $F$-endomorphism of $R$ viewed as
an $F$-module is multiplication by some element $r$ from $R$. By
density of the base algebra $B_\Lambda$ in its completion, every
endomorphism of $R$ is uniquely determined by its action on
$B_\Lambda$. It is therefore enough to show that a given
endomorphism $\varphi$ of $R$ acts as multiplication by some $r
\in R$ on $B_\Lambda$. It is our first aim to show that such
$\varphi$ acts as multiplication on each $x_\alpha $ for $\alpha <
\lambda$. Therefore we need the following

\begin{definition}
Let $W$ be a subset of $\lambda$. We say that $W$ is {\bf closed}
if for every $\alpha \in W$ we have
\[
x_\alpha  \in R_W^{\alpha}:=F[ y_{\eta n x_{\beta}}, z_{\tau} :
 \eta \in V_{\beta}, \tau \in T_{\beta}, \beta \in W, \beta <
 \alpha, n \in \omega].
 \]
 Moreover, we let $R_W:=F[ y_{\eta n x_{\beta}}, z_{\tau} :
 \eta \in V_{\beta}, \tau \in T_{\beta}, \beta \in W, n \in \omega].$
 \end{definition}

We have a first lemma.

\begin{lemma}
\label{closed} Let $W$ be a finite subset of $\lambda$. Then there
exists a finite and closed subset $W^{\prime}$ of $\lambda$
containing $W$.
\end{lemma}

\proof We prove the claim by induction on $\gamma=\max(W)$. If
$\gamma=0$, then $W=\{ 0 \}$, $R_W=F$, $x_0=0$ and there is
nothing to prove. Hence assume that $\gamma>0$. Thus $x_{\gamma}
\in R_{\gamma}=F[y_{\eta n x_{\beta}}, z_{\tau} : n \in \omega,
\eta \in V_{\beta}, \tau \in T_{\beta}, \beta < \gamma]$ and
therefore there exists a finite set $Q \subseteq \gamma$ such that
\[ x_{\gamma} \in F[ y_{\eta n x_{\beta}}, z_{\tau} :
n \in \omega, \eta \in V_{\beta}, \tau \in T_{\beta}, \beta \in
Q].\] Letting $Q_1=Q \cup (W \backslash \{ \gamma \})$ it follows
that $\max(Q_1) < \gamma$. Thus by induction there exists a closed
and finite $Q_2 \subseteq \lambda$ containing $Q_1$. It is now
easy to see that $W^{\prime}=Q_2 \cup \{\gamma \}$ is as
required. \qed

A closed and finite subset $W$ of $\lambda$ gives rise to a nice
presentation of the elements in $R_W$.

\begin{lemma}
\label{closedpresentation} Let $W$ be a closed and finite subset
of $\lambda$ and $r \in R_W$. Then there exists an integer $m_*^r
\in \N$ such that for every $n \geq m_*^r$ we have $r \in
F[y_{\eta n x_{\beta}}, z_{\tau} : \eta \in V_{\beta}, \tau \in
T_{\beta}, \beta \in W ]$.
\end{lemma}

\proof We induct on the cardinality of the finite set $W$. If
$|W|=0$, then $W=\emptyset$ and there is nothing to prove. Thus
assume that $|W| >0$ and let $\gamma$ be maximal in $W$. It is
easy to see that $W^{\prime} = W \backslash \{ \gamma \}$ is
still closed and of course finite. Moreover, by the definition of
closeness we have $x_{\delta} \in R_{W^{\prime}}$ for all $\delta
\in W$. By induction it follows that $x_{\delta} \in F[y_{\eta n
x_{\beta}}, z_{\tau} : \eta \in V_{\beta}, \tau \in T_{\beta},
\beta \in W^{\prime} ]$ for some fixed $m_*^{\delta}$ and every $n
\geq m_*^{\delta}$ ($\delta \in W$). Now let $r \in R_W$. Then $r$
can be written as
\[ r=\sigma(\{ y_{\eta_{r,l} k_{r,l} x_{\beta_{r,l}}},
z_{\tau_{r,j}} : \eta_{r,l} \in V_{\beta_{r,l}}, \tau_{r,j} \in
T_{\beta_{r,j}}, l < l_r, j < j_r \}) \] for some $l_r, j_r \in
\N$, $\beta_{r,l}, \beta_{r,j} \in W$ and $\eta_{r,l} \in
V_{\beta_{r,l}}$, $\tau_{r,j} \in T_{\beta_{r,j}}$. Let
$m_*^r=\max(\{m_*^{\delta}, k_{r,l} : l < l_r, \delta \in W \})$.
Using the (\ref{equat}) it follows now easily that $r \in
F[y_{\eta} n x_{\beta}, z_{\tau} : \eta \in V_{\beta}, \tau \in
T_{\beta}, \beta \in W ]$ for every $n \geq m_*^r$. \qed

We are now ready to show that every endomorphism of $R$ acts as
multiplication on each of the $x_\alpha $'s.

\begin{definition}
Let $R=\bigcup\limits_{\alpha < \lambda}R_\alpha $ be as above.
Then we define \[G_\alpha =\left<  y_{\eta nx_\alpha }, z_{\tau} :
\eta \in V_\alpha , \tau \in T_\alpha , n \in \omega \right>_F,\]
a submodule of the $F$-module $R_\alpha $ for every $\alpha <
\lambda$.
\end{definition}

Clearly it is enough to show that every homomorphism from
$G_\alpha $ to $R^+$ must map $x_\alpha $ to a multiple of itself,
hence is multiplication with some element from $R$ on $x_\alpha $.
Note that $x_\alpha  \in G_\alpha $ by (\ref{equat}).

\begin{proposition}
\label{mainprop} Let $h: G_\alpha  \rightarrow R$ be an $F$-
homomorphism. Then $h(x_\alpha ) \in x_\alpha R$.
\end{proposition}

\proof Let $h: G_\alpha  \rightarrow R$ be an $F$-homomorphism and
assume towards contradiction that $h(x_\alpha ) \not\in x_\alpha
R$. For a subset $V \subseteq V_\alpha $ of cardinality $\lambda$
we define the module
\[
G_V=\left< x_\alpha , y_{\eta n x_\alpha } : \eta \in V, n \in
\omega \right>_* \subseteq G_\alpha
\] Note that
$x_\alpha \in G_V$ implies $\{ z_{\eta\restr_n} : \eta \in
V, n \in \omega \} \subseteq G_V$ by (\ref{equat}).\\
Obviously, the set
\[ H=\{ G_V : V \subseteq V_\alpha , |V|=\lambda \}
\]
is not empty since $G_{V_\alpha }$ belongs to $H$. Let $\beta_*$
be the minimum of the set $\{ \beta < \lambda : \exists G_V \in H
\text{ and } h(G_V) \subseteq R_{\beta} \}$ and choose any $G_V
\in H$ such that $h(G_V) \subseteq R_{\beta_*}$. We will first
show
that $\beta_* \not= \lambda$.\\
Assume towards contradiction that $\beta_*=\lambda$. For $\eta
\in V$ we define
\[ Y_{\eta}=\{ y_{\eta n x_\alpha } : n < \omega
\}\] which is countable. Let $U=U_V \subseteq V$ be a countable
subset as in equation (\ref{representatives}). Moreover, since
$\lambda$ is regular uncountable, hence $\cf(\lambda)=\lambda
> \aleph_0$ we can find $\beta < \lambda$ such that
\[ h(Y_{\rho}) \subseteq R_{\beta} \]
for all $\rho \in U$. Without loss of generality $\beta$ is a
successor ordinal and $h(x_\alpha ) \in R_{\beta}$. We fix $n_*
\in \omega$ and $\eta \in V$. Let $n_* < n \in \omega$ and choose
by equation (\ref{representatives}) $\rho_n \in U$ such that
$\eta\restr_n = \rho_n\restr_n$. We obtain
\begin{equation}
\label{divi}
 y_{\eta n_* x_\alpha } -
y_{\rho_n n_* x_\alpha }
\end{equation}
\[=\sum\limits_{i \geq n_*} \frac{q_i}{q_{n_*}} \left(
z_{\eta\restr_i} \right) + x_\alpha \sum\limits_{i \geq n_*}
\frac{q_i}{q_{n_*}} \eta(i)-\sum\limits_{i \geq n_*}
\frac{q_i}{q_{n_*}} \left( z_{\rho_n\restr_i} \right) - x_\alpha
\sum\limits_{i \geq n_*} \frac{q_i}{q_{n_*}} \rho_n(i)\]
\[=\sum\limits_{i \geq n+1} \frac{q_i}{q_{n_*}} \left(
z_{\eta\restr_i} \right) + x_\alpha \sum\limits_{i \geq n}
\frac{q_i}{q_{n_*}} \eta(i)-\sum\limits_{i \geq n+1}
\frac{q_i}{q_{n_*}} \left( z_{\rho_n\restr_i} \right) - x_\alpha
\sum\limits_{i \geq n} \frac{q_i}{q_{n_*}} \rho_n(i)\] which is
divisible by $s_{n-1}$. Therefore the image under $h$
\[ h(y_{\eta n_* x_\alpha } - y_{\rho_n n_* x_\alpha }) \]
is divisible by all $s_{n-1}$ for $n_* < n < \omega$. Since
$h(y_{\rho_n n_* x_\alpha }) \in R_{\beta}$ by the choice of
$\rho_n \in U$ it follows that \[ h(y_{\eta n_* x_\alpha }) +
R_{\beta} \in R/R_{\beta} \] is divisible by infinitely many
$s_n$. Hence we obtain $h(y_{\eta n_* x_\alpha }) \in R_{\beta}$
since $R/R_{\beta}$ is $\aleph_1$-free by Lemma \ref{aleph1frei}.
Since $n_*$ was chosen arbitrarily we conclude that for all $\eta
\in V$, $h(Y_{\eta}) \subseteq R_{\beta}$ and thus $h(G_V)
\subseteq R_{\beta}$ which contradicts the minimality of
$\beta_*$. Therefore $\beta_* \not=\lambda$. \\
We now work with the elements $y_{\eta o x_\alpha }$ for $\eta \in
V$. Since $h(G_V) \subseteq R_{\beta_*}$ we can find presentations
\[
h(y_{\eta o x_\alpha }) = \sigma_{\eta}( \{ y_{\nu_{\eta,l}
m_{\eta,l} x_{\beta_{\eta,l}}}, z_{\tau_{\eta,k}} : l < l_{\eta},
k< k_{\eta} \})
\]
for every $\eta \in V$ and suitable $\beta_{\eta,l},
\beta_{\eta,k} < \beta_*$, $\nu_{\eta,l} \in V_{\beta_{\eta,l}}$
and $\tau_{\eta,k} \in T_{\beta_{\eta,k}}$. Recall that
$\sigma_{\eta}$ is indexed by $\eta$ since the ``structure'' of
the polynomial depends on $\eta \in V$. For notational simplicity
we will assume that all $\beta_{\eta,l}$ and $\beta_{\eta,k}$ are
pair-wise distinct. By a pigeon hole argument we may assume
without loss of generality that for all $\eta \in V$ we have
$l_{\eta} =l_*$ and $k_{\eta}=k_*$ for some fixed $l_*, k_* \in
\N$. Moreover, since $F$ is countable, we may assume that all
$\sigma_{\eta}$ are independent of $\eta$, say
$\sigma_{\eta}=\sigma$. Hence we get
\[
h(y_{\eta o x_\alpha }) = \sigma( \{ y_{\nu_{\eta,l} m_{\eta,l}
x_{\beta_{\eta,l}}}, z_{\tau_{\eta,k}} : l < l_*, k< k_* \})
\]
 We put
\[ W_{\eta}=\{ \beta_{\eta,l}, \beta_{\eta,k} : l < l_*, k < k_* \} \]
which is a finite subset of $\lambda$ for every $\eta \in V$. By
Lemma \ref{closed} we may assume without loss of generality that
$W_{\eta}$ is already closed $(\eta \in V)$. Moreover, we may
assume that $h(x_\alpha ) \in R_{W_{\eta}}$ for all $\eta \in V$
by possibly enlarging $W_{\eta}$. Since $\beta_*$ is less than
$\lambda$ and $\lambda$ is regular there is a finite and closed
subset $W=\{ \beta_l, \beta_k : l< l_*, k< k_*\}$ of $\lambda$
such that $W_{\eta}=W$ for $\eta \in V^{\prime}$ for some
$V^{\prime} \subseteq V$ of cardinality $\lambda$. Without loss of
generality we will assume $V=V^{\prime}$. Let $m_{\eta} \in \N$
such that $m_{\eta} > lg(\tau_{\eta,k})$ for all $\eta \in V$ and
$k < k_*$. Again, possibly shrinking the set $V$ but preserving
its cardinality we may assume by a pigeon hole argument that
$m_{\eta}=m_1$ for some fixed $m_1 \in \Z$ and all $\eta \in V$.
Now we apply Lemma \ref{closedpresentation} to obtain $h(y_{\eta o
x_\alpha }) \in F[y_{\eta n_{\eta} x_{\beta}}, z_{\tau} : \eta \in
V_{\beta}, \tau \in T_{\beta}, \beta \in W ]$ for $\eta \in V$ and
some $n_{\eta} \in \N$. Once more applying a pigeon hole argument
we may assume that $n_{\eta}=n_*$ for some fixed $n_* \in \N$ and
all $\eta \in V$. We let $m_*=\max\{n_*, m_1\}$ to find new
presentations
\begin{equation}
\label{new} h(y_{\eta o x_\alpha }) = \sigma( \{ y_{\nu_{\eta,l}
m_* x_{\beta_l}}, z_{\tau_{\eta, k}} : l < l_*, k< k_* \})
\end{equation}
for every $\eta \in V$ and $\beta_l, \beta_k \in W$,
$\nu_{\eta,l} \in V_{\beta_l}$ and $\tau_{\eta, k} \in
T_{\beta_k}$. Moreover, we now have \[ lg(\tau_{\eta,k})
\leq m_* \text{ for all } \eta \in V \text{ and } k < k_*.
\]
The reader may notice that when obtaining equation (\ref{new}) the
polynomial $\sigma$ and the natural number $k_*$ may become
dependent on $\eta$ again but a pigeon hole argument allows us to
unify them again and for notational reasons we stick to $\sigma$
and $k_*$. Another pigeon hole argument, using the fact that
$T_\alpha $ is countable,  allows us to assume that
$\tau_{\eta,k}=\tau_k$ for all $\eta \in V$ and $k< k_*$, hence
\[
h(y_{\eta o x_\alpha }) = \sigma( \{ y_{\nu_{\eta,l} m_*
x_{\beta_l}}, z_{\tau_k} : l < l_*, k< k_* \}).
\]
Finally, increasing $m_*$ (and unifying $\sigma$ and $k_*$ again)
we may assume that all $\nu_{\eta,l}\restr_{m_*}$ are different
$(l < l_*)$ and that
 \begin{equation}
 \label{nottau}
 \nu_{\eta,l}\restr_{m_*} \not= \tau_k
 \end{equation}
  for all $\eta \in
V$ and $l< l_*, k< k_*$. Once more using a pigeon hole argument and
the countability of the trees $T_{\beta_l}$ we may assume that
\begin{equation}
\label{independent} \nu_{\eta,l}\restr_{m_*} = \bar{\tau_l} \in
T_{\beta_l}
\end{equation}
independent of $\eta \in V$ for all $l< l_*$. Hence, $\tau_k \not=
\bar{\tau_l}$ for all $l < l_*$ and $k< k_*$. Last but not least,
since $W$ is closed and $h(x_\alpha ) \in R_W$ we can find
presentations
\[
h(x_{\beta}) = \sigma_{\beta}( \{ y_{\nu_{\beta,l} m_*
x_{\beta_l}}, z_{\tau_{\beta,k}} : l< l_{\beta}, k< k_{\beta} \})
\]
for every $\beta \in W \cup \{ \alpha \}$ and suitable
$l_{\beta}, k_{\beta} \in \N$, $\beta_l, \beta_k \in W$.
Obviously we may assume, once more increasing $m_*$, that
\begin{equation}
\label{notbeta} \nu_{\beta, l}\restr_{m_*} \not=
\nu_{\beta^{\prime}, l^{\prime}}\restr_{m_*} \text{ and }
\nu_{\beta, l}\restr_{m_*} \not= \bar{\tau_j}
\end{equation} for all $\beta, \beta^{\prime} \in W \cup \{
\alpha\}$, $l < l_{\beta},
l^{\prime} < l_{\beta^{\prime}}, j < l_*$.\\
Now choose any $n_* > m_*$ such that
\begin{enumerate}
\item $n_* > \sup(C_{\beta} \cap C_{\beta^{\prime}})$ for all
$\beta \not= \beta^{\prime} \in W \cup \{ \alpha \}$;
\item $s_{n_*}$ is relatively prime to all coefficients in $\sigma$;
\item $s_{n_*}$ is relatively prime to all coefficients in
$\sigma_{\beta}$ for all $\beta \in W \cup \{ \alpha \}$.
\end{enumerate}
Let $n> n_*$ be arbitrary such that for some $\eta_1, \eta_2 \in
V$ we have $\br(\eta_1, \eta_2)=n+1$. Note that by the
uncountability of $V$ there must be infinitely many such $n$, say
for $n \in U \subseteq \omega$, $|U|=\omega$ . An easy calculation
using (\ref{equat}) shows that
\[
y_{\eta_1 o x_\alpha } - y_{\eta_2 o x_\alpha }=
\left(\prod\limits_{l \leq n} s_l \right)(y_{\eta_1 n x_\alpha } -
y_{\eta_2 n x_\alpha })
\]
and hence modulo $s_{n+1}R$ we obtain
\begin{equation}
\label{modulo} y_{\eta_1 o x_\alpha } - y_{\eta_2 o x_\alpha }
\equiv \left( \prod\limits_{l \leq n} s_l \right) x_\alpha  \text{
mod } s_{n+1}R
\end{equation}
since $\br(\eta_1, \eta_2)=n+1$.\\
We now distinguish three cases.\\
{\it Case 1:} Assume that for some $l < l_*$ we have
$\br(\nu_{\eta_1,l}, \nu_{\eta_2,l}) > n+1$. Then clearly
(\ref{equat}) we have
\[ y_{\nu_{\eta_1,l} m_* x_{\beta_l}} - y_{\nu_{\eta_2,l} m_*
x_{\beta_l}} \equiv 0 \text{ mod } s_{n+1}R. \] {\it Case 2:}
Assume that for some $l < l_*$ we have $\br(\nu_{\eta_1,l},
\nu_{\eta_2,l}) = n+1$. Then, again by (\ref{equat}) we have
\[ y_{\nu_{\eta_1,l} m_* x_{\beta_l}} - y_{\nu_{\eta_2,l} m_*
x_{\beta_l}} + s_{n+1}R \in x_{\beta_l}R/s_{n+1}R. \] We will
show that $\beta_l = \alpha$. But this follows from $n> n_* >
\sup(C_{\beta} \cap C_{\beta^{\prime}})$ for all $\beta \not=
\beta^{\prime} \in W \cup \{ \alpha \}$. Hence $n+1$ can not be
the splitting point of pairs of branches from different levels
$\alpha$ and $\beta_l$. Note that also $\br(\eta_1,\eta_2)=n+1$.
Thus $\beta_l =\alpha$ and we obtain
\[ y_{\nu_{\eta_1,l} m_* x_{\beta_l}} - y_{\nu_{\eta_2,l} m_*
x_{\beta_l}} + s_{n+1}R \in x_\alpha R/s_{n+1}R. \] {\it Case 3:}
Assume that for some $l < l_*$ we have $k=\br(\nu_{\eta_1,l},
\nu_{\eta_2,l}) < n+1$, hence $m_* < k$ by equation
(\ref{independent}). Then certainly by (\ref{equat}) and the
choice of $n$ we see that $y_{\nu_{\eta_1,l} n x_\alpha }$ appears
in some monomial in $h(y_{\eta_1 o x_\alpha } - y_{\eta_2 o
x_\alpha })$ with coefficient relatively prime to $s_{n+1}$. By
an easy support argument (restricting to $\nu_{\eta_1,l}\restr_k$
and using  (\ref{nottau}), (\ref{independent}) and
(\ref{notbeta})) this monomial can not appear in $h(x_\alpha )$,
hence
\[ h(y_{\eta_1 o x_\alpha } - y_{\eta_2 o x_\alpha }) - \left(
\prod\limits_{l \leq n} s_l \right) h(x_\alpha ) \not\equiv 0
\text{ mod } s_{n+1}R \] contradicting equation (\ref{modulo}). \\
Therefore, for all $n \in U$ we obtain
\[\left(
\prod\limits_{l \leq n} s_l \right) h(x_\alpha ) \in s_{n+1}R +
x_\alpha R. \] As remarked above the set $U$ is infinite and hence
we obtain
\[ h(x_\alpha ) \in \bigcap\limits_{n \in U} s_{n+1}R + x_\alpha R. \]
We will show that \[ \bigcap\limits_{n \in U} s_{n+1}R + x_\alpha
R =x_\alpha R \] which then implies $h(x_\alpha ) \in x_\alpha R$
and finishes the proof. Recall that $x_\alpha $ is of the form
$x_\alpha =\sum\limits_{m \in M}m$ for some finite subset $M$ of
$M(T_\Lambda)$. We induct on the size $k=|M|$ of $M$. If $k=0$,
then $x_\alpha =0$ and there is nothing to prove. Thus assume that
$k>0$. Let $y \in \bigcap\limits_{n \in U} s_{n+1}R + x_\alpha R$
and choose $f_n, r_n \in R$ for $n \in U$ such that
\begin{equation}
\label{topology} y - s_nf_n = x_\alpha r_n.
\end{equation}
Pick any $m \in M$ and let $M^{\prime}=M \backslash \{ m\}$,
$x_\alpha ^{\prime}=x_\alpha  - m$. We now collect the parts of
$y$ and $f_n$ which have $m$ in their support and therefore write
$y=y_1 + y_2$ and $f_n=f_{n, 1} + f_{n, 2}$ where $m \in [y_1],
[f_{n,1}]$ and $m \not\in [y_2], [f_{n,2}]$ for all $n \in U$.
Restricting equation (\ref{topology}) to $m$ gives
\[ y_1 - s_nf_{n,1}= m r_n \text{ and } y_2 - s_n f_{n,2} =
x_\alpha ^{\prime}r_n \] for all $n \in U$. By the induction
hypothesis we now obtain that there exist $y_1^{\prime},
y_2^{\prime} \in R$ such that
\begin{equation}
\label{primes} y_1 = my_1^{\prime} \text{ and } y_2=x_\alpha
^{\prime}y_2^{\prime}.
\end{equation}
Combining equations (\ref{topology}) and (\ref{primes}) we obtain
\[y - s_nf_n=(y_1 +y_2) -s_nf_n = (my_1^{\prime} +
x_\alpha ^{\prime}y_2^{\prime})- s_nf_n = (m + x_\alpha
^{\prime})r_n = x_\alpha r_n\] for all $n \in U$. Now an easy
support argument shows that $y_1^{\prime}=y_2^{\prime}=y^{\prime}$
for some fixed $y^{\prime} \in R$ and therefore $y=(m + x_\alpha
^{\prime})=x_\alpha y^{\prime} \in x_\alpha R$.

 \qed

We are now ready to prove that $R$ is an $E(F)$-algebra.

\begin{mtheorem}
Let $F$ be a countable principal ideal $p$-domain with identity,
not a field, and let $\aleph_1 \leq \lambda \leq 2^{\aleph_0}$ be
a regular cardinal. If $R=\bigcup\limits_{\alpha < \lambda}
R_\alpha $ is the $F$-algebra constructed above, then $R$ is an
$\aleph_1$-free $E(F)$-algebra of cardinality $\lambda$.
\end{mtheorem}
\proof Let $h$ be any $F$-endomorphism of $R$ viewed as
$F$-module. We have to show that $h$ acts as multiplication by
some element $b \in R$. Restricting $h$ to $G_\alpha $ and
applying Proposition \ref{mainprop} we know that for every $\alpha
< \lambda$ there exists an element $b_\alpha  \in R$ such that
\[ h(x_\alpha )=x_\alpha b_\alpha . \]
Since the elements $x_\alpha $ ($\alpha < \lambda$) run through
all pure elements of $B_{\Lambda}$ we therefore obtain that for
any $r \in B_{\Lambda}$ there is some $b_r \in R$ such that
\[ h(r)=rb_r. \]
Note that for every element $r\in B_{\Lambda}$ there is $\alpha <
\lambda$ and $f \in F$ such that $r=x_\alpha f$ and hence
$h(r)=h(x_\alpha f)=x_\alpha b_\alpha f=rb_\alpha $.\\ Now let
$U_\alpha $ be a countable subset of $V_\alpha $ for every
$\alpha < \lambda$ as in equation (\ref{representatives}). Let \[
R_\alpha ^*=F[y_{\eta n x_{\beta}}, z_{\tau} : \eta \in
U_{\beta}, \tau \in T_{\beta}, \beta < \alpha, n \in \omega] \] a
countable subalgebra of $R_\alpha $. Since $\lambda$ is regular
uncountable there exists for every $\alpha < \lambda$ an ordinal
$\gamma_\alpha  < \lambda$ such that
\[ h(R_\alpha ^*) \subseteq R_{\gamma_\alpha }. \]
We put $C=\{ \delta < \lambda : \forall (\alpha < \delta)
(\gamma_\alpha  < \delta) \}$ which is a closed unbounded subset
(cub) of $\lambda$. Without loss of generality $C$ consists of
limit ordinals. Let $\delta \in C$, then similar arguments as in
the proof of Proposition \ref{mainprop} after equation
(\ref{representatives}), using the fact that $R/R_{\delta}$ is
$\aleph_1$-free show that $h(R_{\beta}) \subseteq R_{\delta}$ for
every $\beta < \delta$ and hence \[ h(R_{\delta}) \subseteq
R_{\delta}. \] Let us assume for the moment that there is some
$\delta_* \in C$ such that for every $r \in B_{\Lambda}$ we have
$b_r \in R_{\delta_*}$. Choose $r_1 \not= r_2 \in B_{\Lambda}$ and
assume that $b_{r_1} \not = b_{r_2}$. Choose $\delta_* < \delta
\in C$ such that $r_1, r_2 \in R_{\delta}$. Let $\tau \in
T_{\delta}$ be arbitrary such that $\tau \not\in ([r_1] \cup
[r_2]$). Then
\begin{equation}
\label{finalequ} b_{\tau}\tau+b_{r_1}r_1 = h(\tau) + h(r_1) =
h(\tau + r_1) = b_{\tau+r_1} (\tau+r_1) = b_{\tau+r_1}\tau +
b_{\tau+r_1}r_1.
\end{equation}
Now note that $R_{\delta}$ is an $R_{\delta_*}$-module and that
$R/R_{\delta}$ is torsion-free as an $R_{\delta_*}$-module.
Moreover, $b_\tau, b_{r_1}$ and $b_{\tau + r_1}$ are elements of
$R_{\delta_*}$, hence $\tau$ is not in the support of either of
them. Thus restricting equation (\ref{finalequ}) to $\tau$ we
obtain
\[ b_\tau \tau = b_{\tau + r_1} \tau \]
and therefore $b_\tau = b_{\tau + r_1}$. Now equation
(\ref{finalequ}) reduces to $b_{r_1}r_1 = b_{\tau + r_1}r_1$ and
since $R$ is a domain we conclude $b_{r_1}=b_{\tau +r_1}$. Hence
$b_{r_1}=b_\tau$. Similarly we obtain $b_{r_2}=b_\tau$ and
therefore $b_{r_1}=b_{r_2}$ which contradicts our assumption. Thus
$b_r=b$ for every $r \in B_{\Lambda}$ and some fixed $b \in R$ and
therefore $h$ acts as multiplication by $b$ on $B_{\Lambda}$ and thus by density also on $R$. \\
It remains to prove that there is $\delta_* < \lambda$ such that
for every $r \in B_{\Lambda}$ we have $b_r \in R_{\delta_*}$.\\
Assume towards contradiction that for every $\delta \in C$ there
is some element $r_{\delta} \in B_{\Lambda}$ such that
$b_{\delta}=b_{r_{\delta}} \not\in R_{\delta}$. We may write
$r_{\delta}$ and also $b_{r_{\delta}}$ as elements in some
polynomial ring over $R_{\delta}$, i.e. we write
\[ r_{\delta} = \sigma_{r_{\delta}} (x_i^{\delta} : i < i_{r_{\delta}} )
\text{ and } b_{\delta} = \sigma_{b_{\delta}}(\tx_i^{\delta}  : i
< i_{b_{\delta}}) \] where $\sigma_{r_{\delta}}$ and
$\sigma_{b_{\delta}}$ are polynomials over $R_{\delta}$ and the
$x_i^{\delta}$'s and $\tx_i^{\delta}$ are independent elements
over $R_{\delta}$. By a pigeon hole argument we may assume that
for all $\delta \in C$ we have $i_{r_{\delta}}=i_r$ and
$i_{b_{\delta}}=i_b$ for some fixed $i_r, i_b \in \N$. Now choose
$n < \omega$ and note that $\varphi: \bigcup\limits_{\alpha <
\lambda}R_\alpha /s_nR_\alpha  \rightarrow \bigcup\limits_{\alpha
< \lambda} (R_\alpha ^*+s_nR)/s_nR$ is an epimorphism. Let
$\bar{\sigma}_{r_{\delta}}$ and $\bar{\sigma}_{b_{\delta}}$ denote
the images of the polynomials $\sigma_{r_{\delta}}$ and
$\sigma_{b_{\delta}}$ under $\varphi$. Since
$\bigcup\limits_{\alpha < \delta} (R_\alpha ^*+s_nR)/s_nR$ is of
cardinality less than $\lambda$ (for every $\delta < \lambda$) and
$C$ consists of limit ordinals an easy application of Fodor's
lemma shows that the mapping $\phi: C \rightarrow R/s_nR$, $\delta
\mapsto (\bar{\sigma}_{r_{\delta}}, \bar{\sigma}_{b_{\delta}})$ is
constant on some stationary subset $C^{\prime}$ of $C$ and without
loss of generality we may assume that $C=C^{\prime}$.\\
Now fix $\delta \in C$ and choose $\delta_1 \in C$ such that
$x_i^{\delta}, \tx_j^{\delta} \in R_{\delta_1}$ for all $i < i_r,
j< i_b$. Let $\delta_1 < \delta_2 \in C$ and let $R^{\prime}$ be
the smallest polynomial ring over $R_{\delta}$ generated by at
least the elements $x_i^{\delta_1}, x_i^{\delta_2}$ and
$\tx_i^{\delta_1}, \tx_i^{\delta_2}$ such that
\[ a_1a_2=a_3 \text{ and } a_2,a_3 \in R^{\prime} \text{ implies }
a_1 \in R^{\prime}. \] Thus $R^{\prime}=R_{\delta}[H]$ as a
polynomial ring where $H \subseteq R \backslash R_{\delta}$
contains the set $\{ x_i^{\delta_1}, x_i^{\delta_2},
\tx_j^{\delta_1}, \tx_j^{\delta_2} : i < i_r, j < i_b \}$. We now
consider the following equation:
\begin{equation}
\label{b-equation} b_{r_{\delta} + r_{\delta_2}}(r_{\delta} +
r_{\delta_2}) = h(r_{\delta} + r_{\delta_2}) =h(r_{\delta}) +
h(r_{\delta_2}) = b_{\delta}r_{\delta} + b_{\delta_2}r_{\delta_2}.
\end{equation}
By the choice of $R^{\prime}$ we get that $b_{r_{\delta} +
r_{\delta_2}} \in R^{\prime}$ since $r_{\delta}, r_{\delta_2},
b_{\delta}, b_{\delta_2} \in R^{\prime}$. Assume that some
$x_i^{\delta}$ appears in the support of $b_{r_{\delta} +
r_{\delta_2}}$, then the product $x_i^{\delta}x_j^{\delta_2}$
appears on the left side (for some $j < i_b$) of the equation
(\ref{b-equation}) but not on the right side - a contradiction.
Similarly, no $x_i^{\delta_2}$ can appear in the support of
$b_{r_{\delta} + r_{\delta_2}}$. Thus we obtain
\[ (b_{r_{\delta} +r_{\delta_2}}  - b_{\delta})r_{\delta} = -(b_{r_{\delta} +
r_{\delta_2}} - b_{\delta_2})r_{\delta_2} \] and therefore
\[ b_{r_{\delta} + r_{\delta_2}} =  b_{\delta} =b_{\delta_2} .\]
Hence $b_{\delta_2} \in R_{\delta_2}$. But this contradicts the
choice of $r_{\delta_2}$. \qed

\begin{corollary}
There exists an almost-free $E$-ring of cardinality $\aleph_1$.
\end{corollary}

\begin{remark}
The authors would like to remark that the Main Theorem could also
be proved for cardinals $\aleph_1 \leq \lambda \leq 2^{\aleph_0}$
which are not regular. But the proof would be much more technical
and complicated, in particular in the case when
$\cf(\lambda)=\omega$. We therefore restricted ourselves to
regular cardinals since our main interest lies in
$\lambda=\aleph_1$.
\end{remark}


 \end{document}